\documentclass[onecolumn]{elsarticle}
\makeatletter
\def\ps@pprintTitle{%
	\let\@oddhead\@empty
	\let\@evenhead\@empty
	\def\@oddfoot{\centerline{\thepage}}%
	\let\@evenfoot\@oddfoot}
\makeatother
\usepackage[english]{babel}
\usepackage{amsfonts}
\usepackage{amsmath}
\usepackage{amscd}
\usepackage{dirtytalk}
\usepackage[active]{srcltx} % SRC Specials: DVI [Inverse] Search
\usepackage{latexsym}
\usepackage{amssymb}
\usepackage[latin1]{inputenc}
\usepackage{enumerate}

\usepackage{pdfsync}
\usepackage{amsthm}
\usepackage{amsmath}
\usepackage{mathtools}
\usepackage{amscd}
\usepackage{latexsym}
\usepackage{amssymb}
\usepackage{enumitem}
\newtheorem{thm}{Theorem}[section]
\newtheorem{prop}[thm]{Proposition}
\newtheorem{cor}[thm]{Corollary}
\newtheorem*{cor*}{Corollary}
\newtheorem{lema}[thm]{Lemma}
\newtheorem*{lema*}{Lemma}

\numberwithin{equation}{section}
\theoremstyle{definition}
\newtheorem*{Def}{Definition}

\newenvironment{dem}{\vspace{1ex}\noindent{\it Proof.}\hspace{0.5em}}
{\hfill\qed\vspace{1ex}}

\newtheorem{example}{Example}
\newtheorem*{obs}{Remark}

\newtheorem*{obs*}{Remark}
\newtheorem*{thm*}{Theorem}
\newtheorem*{prop*}{Proposition}

\newtheoremstyle{dotless}{}{}{}{}{}{}{ }{}
\theoremstyle{dotless}

\newcommand{\PI}[2]{\left\langle \,#1 , #2\, \right\rangle}

\newcommand{\NC}[1]{\Vert #1 \Vert}

\newcommand{\WS}{W_{/ [\St]}}

\newcommand{\K}[2]{[ \,#1 , #2\, ]}

\newcommand{\ra}{\rightarrow}

\newcommand{\St}{\mathcal{S}}

\newcommand{\HH}{\mathcal{H}}

\newcommand{\M}{\mathcal{M}}
\newcommand{\N}{\mathcal{N}}

\newcommand{\Q}{\mathcal{Q}}
\newcommand{\T}{\mathcal{T}}

\newcommand{\mc}[1]{\mathcal{#1}}

\newcommand{\ol}{\overline}

\newcommand{\perpi}{[\perp]}

 \DeclareMathOperator{\tr}{tr}

\begin{document}

\begin{frontmatter}
	
	\title{Weighted operator least squares problems and the $J$-trace in Krein spaces}

	\author[FI,IAM]{Maximiliano Contino\corref{ca}}
	\ead{mcontino@fi.uba.ar}
	
	\author[FI,IAM]{Alejandra Maestripieri}
	\ead{amaestri@fi.uba.ar}
	
		\author[IAM,IV,UNGS]{Stefania Marcantognini}
	\ead{smarcantognini@ungs.edu.ar}

	\cortext[ca]{Corresponding author}
	\address[FI]{%
		Facultad de Ingenier\'{\i}a, Universidad de Buenos Aires\\
		Paseo Col\'on 850 \\
		(1063) Buenos Aires,
		Argentina 
	}
	
	\address[IAM]{%
		Instituto Argentino de Matem\'atica ``Alberto P. Calder\'on'' \\ CONICET\\
		Saavedra 15, Piso 3\\
		(1083) Buenos Aires, 
		Argentina }
	
	\address[IV]{%
		Departamento de Matem\'atica -- Instituto Venezolano de Investigaciones Cient\'ificas \\ Km 11 Carretera Panamericana \\ Caracas, Venezuela
	}

	\address[UNGS]{%
 Universidad Nacional de General Sarmiento -- Instituto de Ciencias \\ Juan Mar\'ia Gutierrez \\ (1613) Los Polvorines, Pcia. de Buenos Aires, Argentina
}
	
\begin{abstract}
	Given $B, C$ and $W$ operators in the algebra $L(\HH)$ of bounded linear operators on the Krein space $\HH,$ the minimization problem 
	$\min \ (BX - C)^{\#}W(BX -C),$ for $X\in L(\HH),$ is studied when the weight $W$ is  selfadjoint. The analogous maximization and min-max problems are also considered. Complete answers to these problems and to those naturally associated to trace clase operators on Krein spaces are given.
	
\end{abstract}

	
\end{frontmatter}

\section{Introduction}

In estimation theory one would like to approximate the values of certain quantities that are not directly observable from the values of some sampled measurements. The solution to the problem of estimating the unobservable quantities given the observable ones depends on the model one uses to describe the relation between them and the optimality criterion one chooses to determine the desired estimates. The weighted least squares method is the standard approach in situations when it may not be feasible to assume that every observation should be treated equally. It works by incorporating a  {\emph{weight}}  to each data point as a way to describe its influence over the estimates. 

The Krein space estimation theory developed by Hassibi et al.  \cite{HassibipartI} has brought into play {\emph{indefinite}} weighted least squares problems. Some of those problems were studied in their ``pointwise'' form, for linear operators on infinite-dimensional spaces in  \cite{GiribetKrein} and, for matrices with complex entries  in \cite{HassibipartII, Hassibietal}. Roughly speaking, if one is given an infinite or finite-dimensional linear space $\HH$, a weight $W,$ bounded linear operators or matrices $B, C$, and a vector $y\in \HH$, then the problem is to find an ``extremal''  vector $x_0\in \HH$ for the quadratic form $[W(Bx - Cy), Bx-Cy]$ with $\K{ \ }{ \ }$ a Krein space inner product on $\HH$. If $R(B)$, the range of $B$, is closed and $W$-nonnegative, the vector $x_0$ one seeks minimizing the above quadratic form is called a {\emph{weighted indefinite least squares solution of}} $Bx = Cy$. 

In this work we look instead for a ``global" solution of the problem, meaning a bounded linear operator $X_0$ acting as a $W${\emph{-inverse of}} $B$. Broadly speaking, we consider a Krein space 
$(\HH, \K{ \ }{ \ })$, a selfadjoint operator $W$ on $\HH$ and bounded linear operators $B, C$ on $\HH$. We then determine whether there exists $X_0$ such that, for each $y \in \HH$, $X_0 y$ is a weighted indefinite least squares solution of $Bx=Cy$.  
For a positive weight $W$, the notion of $W$-inverse was introduced by  Mitra and Rao in the case of matrices \cite{Mitra}, and later on extended to Hilbert space operators in \cite{WGI, Contino}. Here we say that $X_0$ is an {\emph{indefinite minimum solution of $BX-C=0$ with weight $W$}} if $X_0$ realizes the minimum of $(BX-C)^{\#}W(BX-C)$ as $X$ runs over $L(\HH)$, the space of the bounded linear operators on $\HH$, where the order is induced by the cone of $\K{ \ }{ \ }$-positive operators of $L(\HH)$. Necessary and sufficient conditions for the existence of such a solution are given and we show that the solution of $BX - I=0$, if it exists, is none other than the Schur complement of $W$ to $R(B)$; i.e., 
$$W_{/ [R(B)]} = \underset{X \in L(\HH)}{\min} \ (BX-I)^{\#}W(BX-I).$$

Given the $W$-indefiniteness of the range of $B$, it is natural to consider min-max problems. In fact, any factorization of $B$ as the sum of two operators, one with $W$-nonnegative range and the other with $W$-nonpositive range, yields a min-max problem. 
As with the minimization problem, we give necessary and sufficient conditions for the solvability of the min-max problem and we obtain another characterization of the Schur complement. Furthermore, even though the decomposition of $B$
depends on the chosen signature operator $J$, the solutions to the min-max problem does not. 

In the Hilbert space setting an associated minimizing problem can be considered in the context of unitarily invariant norms, particularly, in the $p$-Schatten class norms 
$\| \ \|_p,$ in which case -- and under the assumption that $W$ is positive -- it takes the form of the Procrustes problem $\underset{X \in L(\HH)}{\min} \ \|W^{1/2}(BX-C)\|_p$.  Indeed, these two kinds of problems are closely related, as \cite{Nashed, Gold1, Contino} have shown. Inspired by the work of Kintzel on an indefinite Procrustes problem expressed as a max-min problem on traces of matrices \cite{Ulric}, we define a $J$-trace, $\tr_J$, and study the corresponding min-max problem. We find that, if the problem is solvable for every $C$, the solution is unique and equals  $\tr_J(C^{\#}W_{/[ R(B)]}C)$. In addition, if $\tr_J(T)< \infty$ for some signature operator $J$, then 
$\tr_{J'}(T) <\infty$ for any other signature operator $J',$ though it may happen  that $\tr_J(T) \ne \tr_{J'}(T)$. Consequently, the min-max value for  the $\tr_J$ depends on $J$,  but the set of solutions where this value is attained for each $J,$ is independent of $J$.

The paper may be thought of as the second part of \cite{Contino3}, for it contains the weighted versions of the operator least squares problems we studied there. There the fundamental tool for solving the least squares problems was given by the {\emph{indefinite inverse}}. In this work the {\emph{Schur complement}}, as defined and studied in \cite{Contino4},  plays this role.

The paper has four additional sections. Section 2 fixes notation and recalls the basics of Krein spaces, Section 3 gives a brief account of the fundamental results on the Schur complement from \cite{Contino4}. In Section 4  we turn to weighted least squares problems. Subsection 4.1 is entirely devoted to the weighted min-max problems and contains the main results. Section 5 extends the notion of the trace of an operator to the Krein space setting, and applies the results obtained in the previous section to trace-type min and min-max problems for operators.   

\section{Preliminaries}

We assume that all Hilbert spaces are complex and separable. If $\HH$ is a Hilbert space, $L(\HH)$ stands for the algebra of bounded linear operators on $\HH$ and $L(\HH)^+$ for the cone of positive semidefinite operators in  $L(\HH).$  We write $CR(\HH)$ to indicate the subset of $L(\HH)$ of operators with closed range.

The range and nullspace of any $A \in L(\HH)$ are denoted by $R(A)$ and $N(A)$, respectively. Given a subset $\T \subseteq \HH,$ the preimage of $\T$ under $A$ is denoted by $A^{-1}(\T)$ so $A^{-1}(\T)=\{ h \in \HH: \ Ah \in \T \}.$ Given two operators $S, T \in L(\HH),$ the notation  $T \leq_{\HH} S$ signifies that $S-T \in L(\HH)^+.$ 
For any $T \in L(\HH),$ $\vert T \vert := (T^*T)^{1/2}$ is the modulus of $T$ and $T=U\vert T\vert$ is the polar decomposition of  $T,$ with $U$ the partial isometry such that $N(U)=N(T).$

The direct sum of two closed subspaces $\M$ and $\N$ of $\HH$ is represented by $\M \dot{+} \N.$ 
If $\HH$ is decomposed as $\HH=\M \dot{+} \N,$ the projection onto $\M$ with nullspace $\N$ is denoted by $P_{\M {\mathbin{\!/\mkern-3mu/\!}} \N}$ and abbreviated $P_{\M}$ when $\N = \M^{\perp}.$ $\Q$ indicates  the subset of oblique projections in  $L(\HH),$ namely, $\Q:=\{Q \in L(\HH): Q^{2}=Q\}.$

\subsection*{\textbf{Krein Spaces}}

A linear space $\mathcal H$ endowed with an indefinite inner product (a Hermitian sesquilinear form) $\K{ \ }{ \ }$ is a {\emph{Krein space}} if  
$\mathcal {H}$ is the algebraic direct sum of two subspaces $\mathcal {H}_+$ and $\mathcal {H}_-$ such that:
(1)~$\K{x_+}{x_-} =0$ for every $x_\pm \in \HH_\pm$,
and (2)~$(\HH_+, \K{ \ }{ \ })$ and $(\HH_-, -\K{ \ }{ \ })$  are Hilbert spaces. We write
\begin{equation} \label{fundamentaldecom}
\HH=\HH_+ \ [\dotplus] \ \HH_-
\end{equation}
to indicate that the Krein space $\HH$ is the $\K{ \ }{ \ }$-orthogonal direct sum of $\HH_+$ and $\HH_-$, and we say that (\ref{fundamentaldecom}) is a 
{\emph{fundamental decomposition}} of $\HH$. 

In general, all geometrical notions on a Krein space are to be understood with respect to the indefinite inner product. In particular, the {\emph{orthogonal companion}} of a set $\T$ in  $\HH$, which we denote by $\T^{\perpi}$, is the subspace of those  $h \in \HH$ such that $[h,x] = 0$ for all $x\in \T$.

Every fundamental decomposition $\HH=\HH_+ \ [\dotplus] \ \HH_-$ of a given Krein space $(\HH, \K{ \ }{ \ })$ induces a Hilbert space inner product $\PI{ \ }{ \ }$ on $\mathcal H$. Namely, 
$\langle x , y \rangle :=[x_+,y_+] - [x_-,y_-],$ for $x, y \in \mathcal H$, $x = x_+ + x_-$ and $y = y_++y_-.$ 
In this situation the operator $J$ defined on $x =x_+ +x_-$ by $Jx:= x_+ - x_-$ is called a {\it{signature operator}} of $\mathcal H$. 

If $\HH$ is a Krein space, $L(\HH)$ stands for the vector space of  all the linear operators on $\HH$ which are bounded in an associated Hilbert space $(\HH, \PI{ \ }{ \ }).$ Since the norms generated by different fundamental decompositions of a Krein space $\HH$ are equivalent (see, for instance, \cite[Theorem 7.19]{Azizov}), $L(\HH)$ does not depend on the chosen underlying Hilbert space.

The symbol $T^{\#}$ stands for the $\K{ \ }{ \ }$-adjoint of $T \in L(\HH)$. The set of the operators $T \in L(\HH)$ such that $T=T^{\#}$
is denoted $L(\HH)^s$. If $T\in L(\HH)^s$ and $\K{Tx}{x} \geq 0 \mbox{ for every } x \in \HH,$ $T$ is said to be {\emph{positive}}; the notation  $S \leq T$ signifies that $T-S$ is positive.

Given $W \in L(\HH)^s$ and $\St$ a closed subspace of $\HH,$ we say that $\St$ is $W$-\emph{positive} if 
$\K{Ws}{s} > 0$ for every $s \in \St, \ s\not =0.$ $W$-\emph{nonnegative}, $W$-\emph{neutral}, $W$-\emph{negative} and $W$-\emph{nonpositive} subspaces are defined likewise. If $\St$ and $\T$ are two closed subspaces of $\HH,$ the notation $\St \ [\dotplus]_{W} \  \T$ is used to indicate the direct sum of $\St$ and $\T$ when, additionally, $\K{Ws}{t}=0 \mbox{ for every } s \in \St \mbox{ and } t \in \T.$ 

Standard references on Krein space theory are \cite{AndoLibro}, \cite{Azizov} and \cite{Bognar}. We also refer to \cite{DR} and \cite{DR1} as authoritative accounts of the subject.

\section{Schur complement in Krein Spaces}

In this section we include several results on the Schur complement in Krein spaces that will be useful along the paper. For the proofs the reader is referred to \cite{Contino4}.

The notion of Schur complement (or shorted operator) of $A$ to $\St$ for a positive operator $A$ on a Hilbert space $\HH$ and $\St \subseteq \HH$ a closed subspace, was introduced by M.G.~Krein \cite{Krein}. He proved that the set $\{ X \in L(\HH): \ 0\leq_{\HH} X\leq_{\HH} A \mbox{ and } R(X)\subseteq \St^{\perp}\}$ has a maximum element, which he defined as the {{Schur complement}} $A_{/ \St}$ of $A$ to $\St.$ This notion was later rediscovered by Anderson and Trapp \cite{Shorted2}. 
If $A$ is represented as the $2\times 2$ block matrix $\begin{pmatrix} a & b\\ b^* & c \end{pmatrix}$ with respect to the decomposition of $\HH = \St \oplus \St^{\perp},$ they established the formula 
$$A_{/ \St}= \begin{pmatrix} 0 & 0\\ 0& c - y^*y\end{pmatrix}$$ where $y$ is the unique solution of the equation $b = a^{1/2} x$ such that the range inclusion $R(y) \subseteq \ol{R(a)}$ holds. The solution always exists because $A$ is positive: in this case, $a$ is also positive and the range inclusion $R(b) \subseteq R(a^{1/2})$ holds.

In \cite{AntCorSto06} Antezana et al., extended the notion of Schur complement to any bounded operator $A$ satisfying a weak complementability condition with respect to a given pair of closed subspaces $\St$ and $\T,$ by giving an Anderson-Trapp type formula. In particular, if $A$ is a bounded selfadjoint operator, $\St=\T$ and $A=\begin{pmatrix} a & b\\ b^* & c \end{pmatrix},$ this condition reads $R(b) \subseteq R(\vert a \vert^{1/2}),$ which as noted, is automatic for positive operators. In this case, let $f$ be the unique solution of the equation $b = \vert a \vert^{1/2} x$ such that the range inclusion $R(f) \subseteq \ol{R(a)}$ holds and $a=u\vert a \vert$ the polar decomposition of $a.$ Then, the Schur complement  of $A$ to $\St$ is defined as 
$$A_{/ \St}= \begin{pmatrix} 0 & 0\\ 0& c - f^*uf\end{pmatrix}.$$ 

In \cite{Contino4}, the notions of $\St$-complementability, $\St$-weak complementability and the Schur complement were extended to the Krein space setting in the following fashion.

\begin{Def}  Let $W \in L(\HH)^s$ and $\St$ be a closed subspace of $\HH.$ The operator $W$ is called $\St$-\emph{complementable} if 
	$$\HH=\St + W^{-1}(\St^{\perpi}).$$
\end{Def}

If $W$ is $\St$-complementable then, for any fundamental decomposition $\HH=\HH_+ \ [\dotplus] \ \HH_-$ with signature operator $J,$ we get that $\HH=\St + (JW)^{-1}(\St^{\perp}).$ Therefore, $W$ is $\St$-complementable if and only if the pair $(JW,\St)$ is \emph{compatible} in (the Hilbert space) $(\HH, \PI{ \ }{ \ })$ for any (and then for every) signature operator $J,$ meaning that there exists a projection $Q$ onto $\St,$ such that  $JWQ=Q^{*}JW,$ see \cite{CMSSzeged}.
From this, it follows that $W$ is $\St$-complementable if and only if there exists a projection $Q$ onto $\St$ such that $WQ=Q^{\#}W.$

In a similar way the $\St$-weak complementability in Krein spaces, with respect to a fixed signature operator $J,$ is defined.

\begin{Def} Let $W \in L(\HH)^s$ and $\St$ be a closed subspace of $\HH.$ The operator $W$ is $\St$-\emph{weakly complementable} with respect to a signature operator $J$ if $JW$ is $\St$-weakly complementable in $(\HH, \PI{ \ }{ \ }).$ 
\end{Def}

In this case, if the matrix representation of $JW$ induced by $\St$ is
\begin{equation} \label{Wdes} 
JW=\begin{bmatrix} 
a & b \\ 
b^* & c \\  
\end{bmatrix}, \end{equation} the $\St$-weak complementability of $W$ is equivalent to $R(b)\subseteq R(\vert a\vert^{1/2}).$ 
The $\St$-weak complementability of $W$ does not depend on the signature operator, see \cite[Theorem 4.4]{Contino4}.
Then, we simply say that $W$ is $\St$-weakly complementable, whenever $W$ is $\St$-weakly complementable with respect to a signature operator $J.$

Let $W \in L(\HH)^ s$ and $\St$ a closed subspace of $\HH.$ Then, by applying the spectral theorem for Hilbert space selfadjoint operators to $A=JW,$ with $J$ any signature operator, $\St$ can be decomposed as 
\begin{equation} \label{WdecompSKrein}
\St=\St_+ \ [\dotplus]_{W} \ \St_-,
\end{equation}
where $\St_{+}$ and $\St_-$ are closed, $\St_+$ is $W$-nonnegative, $\St_-$ is $W$-nonpositive and $\St_{+} \perp \St_{-}.$ Notice that the decomposition in \eqref{WdecompSKrein} need not be unique.

The following is a characterization of the $\St$-weak complementability \cite[Proposition 4.7]{Contino4}.
\begin{prop} \label{PropWC} Let $W \in L(\HH)^s$  and $\St$ be a closed subspace of $\HH.$ 
	Suppose that $\St=\St_+ \ [\dotplus]_{W} \ \St_-$ is any decomposition as in \eqref{WdecompSKrein} for some signature operator $J.$ Then the following statements are equivalent:
	\begin{enumerate}
		\item[i)] $W$ is $\St$-weakly complementable,
		\item [ii)] there exist $W_1, W_2, W_3 \in L(\HH)^s,$ $W_2, W_3 \geq 0$ such that $W=W_1+W_2-W_3$ and $\St \subseteq N(W_1),$ $\St_- \subseteq N(W_2),$ $\St_+ \subseteq N(W_3),$  
		\item [iii)] $W$ is $\St_{\pm}$-weakly complementable.
	\end{enumerate}
\end{prop}

\begin{Def} Let $W \in L(\HH)^s,$ $\St$ be a closed subspace of $\HH$ and $J$ a signature operator. Suppose that $W$ is $\St$-weakly complementable. The \emph{Schur complement} of $W$ to $\St$ corresponding to $J$ is
	$$W_{/ [\St]}^J =J (JW)_{ / \St},$$
	and the $\St$-\emph{compression} of $W$ is $W_{ [\St]}^J = W- W_{/ [\St]}^J.$
\end{Def}

In \cite[Theorem 4.5]{Contino4} it was proved that the Schur complement does not depend on the fundamental decomposition of $\HH.$ Henceforth we write $W_{/ [\St]}$ for this operator and $W_{ [\St]}$ for the $\St$-compression. Also, suppose that $\St=\St_+ \ [\dotplus]_{W} \ \St_-$ is any decomposition as in \eqref{WdecompSKrein} for some signature operator $J.$ If $W$ is $\St$-weakly complementable then
\begin{equation} \label{ShortedsupinfKrein}
\WS= (W_{/ [ \St_+]} )_{/ [\St_-]}=(W_{/[ \St_-]})_{/ [\St_+]}.
\end{equation}
Also, if $W=W_1+W_2-W_3$ as in Proposition \ref{PropWC} then 
\begin{equation} \label{ShortedKrein}
\WS= W_1+{W_2}_{/ [ \St_+]}-{W_3}_{/ [ \St_-]}.
\end{equation}
Moreover, if $W$ is $\St$-complementable then 
\begin{equation}\label{ShortedmaxminKrein}
W_{/ [\St]}=W(I-Q), 
\end{equation}
for any projection $Q$ onto $\St$ such that $WQ=Q^{\#}W.$ 

The following result was proved in \cite[Corollary 4.12]{Contino4}.

\begin{prop} \label{ShortedC2} Let $W \in L(\HH)^s$ and $\St$ be a closed subspace of $\HH.$  Suppose that $\St$ is $W$-nonnegative. Then $W$ is $\St$-weakly complementable if and only if there exists $\inf \ \{ E^{\#}WE: E=E^2, \ N(E)=\St\}.$ In this case, $$\WS=\inf \ \{ E^{\#}WE: E=E^2, \ N(E)=\St\}.$$ 
\end{prop}

\section{Weighted least squares problems in Krein spaces}

Consider the following problem:
given the operators $W \in L(\HH)^s,$ $B\in CR(\HH)$ and $C\in L(\HH),$ determine the existence of
\begin{equation}
\underset{X \in L(\HH)}{\min} (BX-C)^{\#}W(BX-C). \label{eq61}
\end{equation}

\begin{Def} Let $W \in L(\HH)^s,$ $B\in CR(\HH)$ and $C\in L(\HH).$ An operator $X_0 \in L(\HH)$ is an {\emph{indefinite minimum solution of $BX-C=0$ with weight $W$}} ($W$-ImS)  if $X_0$ is a solution of Problem \eqref{eq61}. 
\end{Def}

In a similar fashion, the analogous maximization problem can be considered. 
Along this section all the results are stated for problem \eqref{eq61} but similar results hold for the maximum problem.  

\bigskip
Consider $W \in L(\HH)^s,$ $B \in CR(\HH)$ and $C\in L(\HH)$ and define 
\begin{equation} \label{FX}
F(X):=(BX-C)^{\#}W(BX-C).
\end{equation}
%The following  results concerns  \eqref{eq61} 
We begin by giving conditions for the existence of the infimum in $L(\HH)$ of the family $\{ F(X): X \in L(\HH)\}$ when $C=I.$ 

\begin{prop} \label{propinfimum1} Let $W \in L(\HH)^s$ and $B \in CR(\HH)$ such that $R(B)$ is $W$-nonnegative.
	Then the following are equivalent:
	\begin{itemize}
		\item [i)] There exists $\underset{X \in L(\HH)}{\inf} \ (BX-I)^{\#}W(BX-I)=:Z_0 \in L(\HH)$ and $R(B)$ is $Z_0$-nonnegative, 
		\item [ii)] $W$ is $R(B)$-weakly complementable.
	\end{itemize} 
	In this case, $Z_0=W_{/ [R(B)]}.$
\end{prop}
\begin{dem} Suppose that $W$ is $R(B)$-weakly complementable. Let $F(X)$ be as in \eqref{FX} for $C=I.$ Then, for any $X\in L(\HH),$ 	
	$F(X)=W_{/[R(B)]}+(BX-I)^{\#}W_{[R(B)]}(BX-I) \geq W_{/[R(B)]},$ because $R(B) \subseteq N(W_{/[R(B)]})$ and the fact that $R(B)$ is $W$-nonnegative yields $W_{[R(B)]} \geq 0.$  Hence $W_{/[R(B)]}$ is a lower bound of $\{ F(X) :X \in L(\HH) \}.$  Let $T \in L(\HH)$ be any other lower bound of $F(X).$
	
	In particular, given $E \in \Q$ such that $R(I-E)=R(B),$ by Douglas' Lemma \cite{Douglas}, there exists $X_0 \in L(\HH)$ satisfying $I-E=BX_0;$ i.e., such that $-E=BX_0-I.$ Then
	$$T\leq E^{\#}WE  \mbox { for every } E \in \Q \mbox{ such that } N(E)=R(B).$$
	By Proposition \ref{ShortedC2},
	$$T \leq \inf \  \{ E^{\#}WE: \ E \in \Q, \ N(E)=R(B) \}=W_{/ [R(B)]}.$$
	Therefore, $W_{/ [R(B)]}=\underset{X \in L(\HH)} {\inf} \ F(X)$ and, since $R(B) \subseteq N(W_{/[R(B)]}),$ $R(B)$ is $W_{/ [R(B)]}$-nonnegative.
	
	Conversely, if $Z_0$ exists and $R(B)$ is $Z_0$-nonnegative, then taking $X=0,$ the inequality  $Z_0 \leq W$ shows that $Z_0 \in L(\HH)^s.$ 
	As before, $$Z_0\leq E^{\#}WE  \mbox { for every } E \in \Q \mbox{ such that } N(E)=R(B).$$
	
	Fix a signature operator $J$ and let $(\HH, \PI{ \ }{ \ })$ be the corresponding Hilbert space; consider $E=P_{R(B)^{\perp}}.$ Since $Z_0 \in L(\HH)^s,$ $(JZ_0)^*=JZ_0$ and 
	\begin{equation} \label{infR(B)}
	JZ_0\leq_{\HH} P_{R(B)^{\perp}}JWP_{R(B)^{\perp}}.
	\end{equation}
	Let
	$JW=\begin{bmatrix} 
	a & b \\ 
	b^* & c \\  
	\end{bmatrix}$  and 
	$JZ_0=\begin{bmatrix} 
	z_{11} & z_{12} \\ 
	z_{12}^* & z_{22}\\  
	\end{bmatrix}$ be the matrix representation of $JW$ and $JZ_0$ induced by $R(B),$ respectively. By \eqref{infR(B)}, 
	$$P_{R(B)^{\perp}}JWP_{R(B)^{\perp}}-JZ_0=\begin{bmatrix} 
	-z_{11} & -z_{12} \\ 
	-z_{12}^* & c-z_{22}\\  
	\end{bmatrix} \geq_{\HH}0.$$ Then, $z_{11} \leq_{\HH} 0$ and $R(z_{12}) \subseteq R((-z_{11})^{1/2}).$ Since $R(B)$ is $Z_0$-nonnegative, $z_{11} \geq_{\HH}0.$ So $z_{11}=z_{12}=z_{12}^*=0$ and $R(JZ_0) \subseteq R(B)^{\perp}$ or equivalently, $R(Z_0) \subseteq R(B)^{\perpi}.$ Therefore, $W=(W-Z_0)+Z_0,$ with $W-Z_0 \geq 0$ and $R(Z_0) \subseteq R(B)^{\perpi}.$ Then, by Proposition \ref{PropWC}, $W$ is $R(B)$-weakly complementable. 
\end{dem}

\begin{cor} \label{corinfimum1} Let $W \in L(\HH)^s$ and $B \in CR(\HH)$ such that $R(B)$ is $W$-nonnegative and $W$ is $R(B)$-weakly complementable. Then, for every $C \in L(\HH),$
	$$\underset{X \in L(\HH)}{\inf} \ (BX-C)^{\#}W(BX-C) =C^{\#}W_{/ [R(B)]}C.$$
\end{cor}

\begin{dem} If $W \geq 0,$ by \cite[Lemma 4.1]{Contino}, 
	$$\inf \  \{ C^{\#}E^{\#}WEC: \ E \in \Q, \ N(E)=R(B) \}=C^{\#}W_{/ [R(B)]}C.$$
	By Proposition \ref{PropWC}, $W=W_1+W_2,$ with $R(B) \subseteq N(W_1)$ and $W_2 \geq0.$ 
	Then, given $E \in \Q$ such that $N(E)=R(B),$
	$$C^{\#}E^{\#}WEC=C^{\#}W_1C+C^{\#}E^{\#}W_2EC.$$ Hence
	\begin{equation*} \label{infC}
	\begin{split}
	&\inf \  \{ C^{\#}E^{\#}WEC: \ E \in \Q, \ N(E)=R(B) \}=\\
	&\quad\quad=C^{\#}W_1C+\inf \  \{ C^{\#}E^{\#}W_2EC: \ E \in \Q, \ N(E)=R(B) \}
	\\
	&\quad\quad=C^{\#}W_1C+C^{\#}W_2{_{/ [R(B)]}}C=C^{\#}W_{/ [R(B)]}C.
	\end{split}
	\end{equation*}
	Using this equality, the result follows in a similar way as in the first part of the proof of Proposition \ref{propinfimum1}.
\end{dem}	

The next theorem establishes when the infimum in Proposition \ref{propinfimum1} is attained.

\begin{thm} \label{thminimum2} Let $W \in L(\HH)^s$ and $B \in CR(\HH).$ Then the following are equivalent:
	\begin{itemize}
		\item [i)] there exists a $W$-ImS of  $BX-I=0,$ 
		\item [ii)] 	$R(B)$ is $W$-nonnegative  and  $W$ is $R(B)$-complementable,
		\item [iii)] $R(B)$ is $W$-nonnegative  and the normal equation 
		\begin{equation} \label{NEqW2}
		B^{\#}W(BX-I)=0
		\end{equation} 
		admits a solution.
	\end{itemize}
	In this case,
	$$\underset{X \in L(\HH)}{\min} \ (BX-I)^{\#}W(BX-I)=W_{/ [R(B)]}.$$	
\end{thm}

\begin{dem} 
	$i) \Leftrightarrow iii):$ Suppose that $X_0$ is a $W$-ImS of $BX-I=0.$ Then
	$$ \K{W(BX_0-I)x}{(BX_0-I)x} \leq \K{W(BX-I)x}{(BX-I)x}$$
	$ \mbox{for every } x \in \HH \mbox{ and every } X \in L(\HH).$
	Let $z \in \HH$ be arbitrary. Then, for every $x \in \HH \setminus \{0\},$ there exists $X \in L(\HH)$ such that $z=Xx.$ Therefore
	$$ \K{W(BX_0-I)x}{(BX_0-I)x} \leq \K{W(Bz-x)}{Bz-x}$$
	for every $x, z \in \HH.$ 
	Thus, for every $x \in \HH,$ $X_0x$ is a weighted indefinite least squares solution of $Bz=x.$ 
	So, by \cite[Proposition 3.2]{GiribetKrein} (see also \cite[Chapter I, Theorem 8.4]{Bognar}), $R(B)$ is $W$-nonnegative and $X_0x$ is a solution of $B^{\#}W(By-x)=0$ for every $x \in \HH,$ or equivalently, $X_0$ is a solution of  \eqref{NEqW2}.
	
	The converse follows in a similar way, applying again \cite[Proposition 3.2]{GiribetKrein}.
	
	$ii) \Leftrightarrow iii):$ Suppose that $\HH = R(B) + \ W^{-1}(R(B)^{\perpi}),$ then $R(B^{\#}W) \subseteq R(B^{\#}WB).$ Hence, by Douglas' Lemma, the equation $B^{\#}W(BX-I)=0$ admits a solution.
	The converse follows analogously.
	
	In this case, by Proposition \ref{propinfimum1}, $$\underset{X \in L(\HH)}{\min} (BX-I)^{\#}W(BX-I)=W_{/ [R(B)]}.$$
\end{dem}

The next corollaries follow from Theorem \ref{thminimum2}.

\begin{cor} \label{cormin2} Let $W \in L(\HH)^s,$ $B \in CR(\HH)$ and $C \in L(\HH).$ Then the following are equivalent:
	\begin{itemize}
		\item [i)] there exists a $W$-ImS of  $BX-C=0,$ 
		\item [ii)]  $R(B)$ is $W$-nonnegative  and $R(C) \subseteq R(B)+W^{-1}(R(B)^{\perpi}),$
		\item [iii)] $R(B)$ is $W$-nonnegative  and the normal equation 
		\begin{equation} \label{NEq3}
		B^{\#}W(BX-C)=0
		\end{equation} 
		admits a solution.
	\end{itemize}
	In this case, $X_0$ is a $W$-ImS of $BX-C=0$ if and only $X_0$ is a solution of \eqref{NEq3}.
\end{cor}

\begin{dem} This follows in a similar way as in the proof of Theorem \ref{thminimum2} using the fact that $u$ is a weighted indefinite least squares solution of the equation $Bz=Cx$ if and only if $R(B)$ is $W$-nonnegative and $u$ is a solution of $B^{\#}W(By-Cx)=0,$ see \cite[Proposition 3.2]{GiribetKrein}.
\end{dem}

\begin{cor} \label{corWminimum2} Let $W \in L(\HH)^s$ and $B \in CR(\HH).$
	Then there exists a $W$-ImS of $BX-C=0$ for every $C \in L(\HH)$ if and only if $R(B)$ is $W$-nonnegative and $W$ is $R(B)$-complementable.
	In this case, 	
	$$\underset{X \in L(\HH)}{\min} (BX-C)^{\#}W(BX-C)=C^{\#}W_{/ [R(B)]}C.$$	
\end{cor}

\begin{dem} Suppose that there exists a $W$-ImS of $BX-C=0$ for every $C \in L(\HH).$ Then the conclusion follows by applying Theorem \ref{thminimum2} for $C=I.$ 
	
	Conversely, if $W$ is $R(B)$-complementable and $R(B)$ is $W$-nonnegative, then $B^{\#}W(BX-I)=0$ admits a solution. Therefore  $B^{\#}W(BX-C)=0$ admits a solution for every $C \in L(\HH)$ and, by Corollary \ref{cormin2}, there exists a $W$-ImS of $BX-C=0.$ 
	
	In this case, let $Q \in \Q$ be such that $R(Q)=R(B)$ and $WQ=Q^{\#}W.$ Then, by Douglas' Lemma, there exists $X_0\in L(\HH)$ such that $BX_0=QC.$ Therefore $B^{\#}W(BX_0-C)=B^{\#}W(Q-I)C=0,$ because $R(I-Q)=N(Q)\subseteq N(B^{\#}W)$ \cite[Lemma 3.2]{CMSSzeged}. Then, $X_0$ is a $W$-ImS of $BX-C=0.$ 
	Hence, $\underset{X \in L(\HH)}{\min} F(X)=C^{\#}W_{/ [R(B)]}C,$ since $W_{/ [R(B)]}=W(I-Q),$ by \eqref{ShortedmaxminKrein}.
\end{dem}

\subsection{Weighted Min-Max problems}

A necessary condition for the minimization (maximization) problem to be solvable is that the range of the operator $\!B$ is $W$-nonnegative ($\!W$-nonpositive). 
In what follows, we are interested in posing (and solving) a problem similar to the one in \eqref{eq61}, that does not require the range of $B$ to be $W$-definite in order to admit a solution.
To do so, we begin by expressing the range of $B$ as the sum of suitable $W$-definite subspaces.

For a fix signature operator $J,$ the spectral theorem for Hilbert space selfadjoint operators applied to $JW$ gives that $\St:=R(B)$ can be decomposed as 
$\St=\St_+ \ [\dotplus]_{W} \ \St_-$ (compare with \eqref{WdecompSKrein}). If $P_{\pm}=P_{\St_{\pm}}$ and $B_{\pm}=P_{\pm}B$ then $\St_{\pm}=R(B_{\pm})$ and the following result holds.

\begin{lema} \label{lemmadecomp} Let $W \in L(\HH)^s$ and $B \in CR(\HH).$ Then,  given a signature operator $J,$ 
	$B$ can be written as 
	\begin{equation} \label{DescomposicionB}
	B=B_+ + B_-  
	\end{equation}
	with $R(B_+)$ closed and $W$-nonnegative,  $R(B_-)$ closed and $W$-nonpositive, $R(B_+) \perp R(B_-)$ and $R(B)=R(B_+) \ [\dotplus]_{W} \ R(B_-).$
\end{lema}

Fix a descomposition of $R(B)$ as in \eqref{DescomposicionB} and define
$$F_J(X,Y)=(B_+X+B_-Y-C)^{\#}W(B_+X+B_-Y-C).$$
Notice that $F_J(X,X)=F(X).$

Consider the following problem: determine the existence of 
$$\underset{Y \in L(\HH) }{\max} \left( \underset{X \in L(\HH)}{\min} F_J(X,Y) \right).$$

\begin{prop} \label{propsupinfimum1} Let $W \in L(\HH)^s$ and $B \in CR(\HH)$ such that $W$ is $R(B)$-weakly complementable and $B$  is represented as in \eqref{DescomposicionB} for some signature operator $J.$
	Then, for every $C \in L(\HH),$ 
	$$
	\underset{ Y \in L(\HH)}{\sup} \left( \underset{X \in L(\HH)}{\inf}  F_J(X,Y) \right)=\underset{ X \in L(\HH)}{\inf}  \left( \underset{Y \in L(\HH)}{\sup} F_J(X,Y) \right) =C^{\#} W_{/ [R(B)]} C.
	$$
\end{prop}

\begin{dem} Write $W=W_1+W_2-W_3,$ with $R(B) \subseteq N(W_1),$ $R(B_-) \subseteq N(W_2),$ $R(B_+) \subseteq N(W_3)$ and $W_2, W_3 \geq0$ (see Proposition \ref{PropWC}). Then
	$$F_J(X,Y)\!=\!C^{\#}W_1C\!+\!(B_+X\!-\!C)^{\#}W_2(B_+X\!-\!C)\!-\!(B_-Y\!-\!C)^{\#}W_3(B_-Y\!-\!C).
	$$
	By Proposition \ref{PropWC}, $W$ is $R(B_{\pm})$-weakly complementable. Also, $W$ is $R(B_+)$-weakly complementable and $R(B_+)$ is $W$-nonnegative if and only if $W_2$ is $R(B_+)$-weakly complementable and $R(B_+)$ is $W_2$-nonnegative. Applying Corollary \ref{corinfimum1}, $$\underset{X \in L(\HH)}{\inf} (B_+X-C)^{\#}W_2(B_+X-C)=C^{\#}{W_2}_{ / [R(B_+)]}C.$$
	Therefore, for each $Y \in L(\HH),$
	$$\underset{X \in L(\HH)}{\inf}F_J(X,Y)=C^{\#}W_1C+C^{\#}{W_2}_{ / [R(B_+)]}C-(B_-Y-C)^{\#}W_3(B_-Y-C).$$
	In the same way, by applying Corollary \ref{corinfimum1} and \eqref{ShortedKrein}
	
	\noindent $\! \underset{Y \in L(\HH)}{\sup}\left( \underset{X \in L(\HH)}{\inf} F_J(X,Y) \right)=C^{\#}W_1C+C^{\#}{W_2}_{ / [R(B_+)]}C-C^{\#}{W_3}_{ / [R(B_-)]}C=\\=C^{\#} W_{/ [R(B)]} C.$
	
	\noindent The second equality can be proved similarly.
\end{dem}

\begin{Def} Let $W \in L(\HH)^s,$ $B \in CR(\HH)$ and $C \in L(\HH).$ Suppose that $B$ is represented as in \eqref{DescomposicionB} for some signature operator $J.$ An operator $Z \in L(\HH)$ is an {\emph{indefinite min-max solution of $BX-C=0$ with weight $W$}} ($W$-ImMS) (corresponding to the decomposition given by $J$) if 
	\begin{equation} \label{eqMinMax2}
	(BZ-C)^{\#}W(BZ-C)=\underset{Y \in L(\HH) }{\max} \left( \underset{X \in L(\HH)}{\min} F_J(X,Y) \right).
	\end{equation}
\end{Def}

When the weight is the identity, it was proved in \cite[Theorem 5.1 and Corollary 5.2]{Contino3}, that an operator $Z \in L(\HH)$ is an $I$-ImMS of $BX-C=0,$ for some fundamental decomposition of $\HH,$ if and only if $$Z=Z_1+Z_2$$ where  $B^{\#}(BZ_1-C)=0$ and  $(BZ_2)^{\#}BZ_2=0.$ Therefore, an $I$-ImMS of $BX-C=0$ is independent of the selected fundamental decomposition of $\HH.$ 
Also, there exists an $I$-ImMS of $BX-C=0$  if and only if  $R(C) \subseteq R(B) + R(B)^{\perpi}.$
A similar result holds for a general weight:

\begin{thm} \label{TeominmaxW} Let $W \in L(\HH)^s,$ $B \in CR(\HH)$ and $C \in L(\HH).$
	An operator $Z$ is a $W$-ImMs  of $BX-C=0$  for some (and, hence, any) fundamental decomposition of $\HH,$ if and only if $$Z=Z_1+Z_2$$ where $B^{\#}W(BZ_1-C)=0$ and $(BZ_2)^{\#}WBZ_2=0.$
\end{thm}
The proof follows from Corollary \ref{cormin2}, using similar arguments to those found in the proof of \cite[Theorem 5.1]{Contino3}.

\begin{obs} Let $W \in L(\HH)^s,$ $B \in CR(\HH)$ and $C \in L(\HH).$ Suppose that $B$ is represented as in \eqref{DescomposicionB} for some signature operator $J.$ Then
	$$\underset{Y \in L(\HH) }{\max} \left( \underset{X \in L(\HH)}{\min} F_J(X,Y) \right)=\underset{X \in L(\HH)}{\min} \left( \underset{Y \in L(\HH) }{\max}\ F_J(X,Y) \right).$$
\end{obs}

This follows from Theorem \ref{TeominmaxW} and using similar arguments to those found in the proof of \cite[Remark after Theorem 5.1]{Contino3}.

\begin{cor} \label{WCorollaryminmax} Let $W \in L(\HH)^s,$ $B \in CR(\HH)$ and $C \in L(\HH).$ Then, there exists a $W$-ImMS of $BX-C=0$ if and only if $R(C) \subseteq R(B) + W^{-1}(R(B)^{\perpi}).$
\end{cor}

\begin{dem} 
	Suppose that  $Z$ is a $W$-ImMs  of $BX-C=0.$ Then, by Theorem \ref{TeominmaxW}, $Z=Z_1+Z_2$ where $B^{\#}W(BZ_1-C)=0$ and $(BZ_2)^{\#}WBZ_2=0.$
	Therefore $$R(C) \subseteq R(B) + W^{-1}(R(B)^{\perpi}).$$
	Conversely, if $R(C) \subseteq R(B) + W^{-1}(R(B)^{\perpi})$ then $R(B^{\#}WC) \subseteq R(B^{\#}WB).$
	By Douglas's Lemma, there exists a solution of the normal equation $B^{\#}W(BX-C)=0,$ say $Z_1 \in L(\HH).$  Put $Z_2=0$ and apply Theorem \ref{TeominmaxW} to get that $Z_1$ is a $W$-ImMs  of $BX-C=0.$
\end{dem}

\begin{cor} \label{CorminmaxregW} Let $W \in L(\HH)^s$ and $B \in CR(\HH).$  Then, there exists a $W$-ImMS of $BX-C=0$ for every $C \in L(\HH)$ if and only if W is $R(B)$-complementable.
	In this case, for every signature operator $J,$  
	$$\underset{Y \in L(\HH)}{\max} \left(  \underset{X \in L(\HH)}{\min} F_J(X,Y)\right)=C^{\#}W_{/ [R(B)]}C=C^{\#}W(I-Q)C,$$
	where $Q$ is any projection onto $R(B)$ such that $WQ=Q^{\#}W.$
\end{cor}

\begin{dem} If W is $R(B)$-complementable then, for every $C \in L(\HH),$ $R(C) \subseteq R(B) + W^{-1}(R(B)^{\perpi})$ and, by Corollary \ref{cormin2}, there exists a $W$-ImMS of $BX-C=0.$
	
	Conversely, assume that, for every $C \in L(\HH)$ there exists a $W$-ImMS of $BX-C=0.$ Set $C=I$ and apply the corollary once again to get that  W is $R(B)$-complementable as $\HH=R(I) \subseteq R(B)+ W^{-1}(R(B)^{\perpi}).$
	
	In this case, like in the proof of Corollary \ref{corWminimum2}, let $Q \in \Q$ be such that $R(Q)=R(B)$ and $WQ=Q^{\#}W.$ Then, by Douglas' Lemma, there exists $Z_1\in L(\HH)$ such that $BZ_1=QC$ and  $B^{\#}W(BZ_1-C)=0.$ Then, by Theorem \ref{TeominmaxW}, $Z_1$ is a $W$-ImMS of $BX-C=0.$
	Therefore,
	$$\underset{Y \in L(\HH)}{\max} \left( \underset{X \in L(\HH)}{\min} F_J(X,Y)\right) =F(Z_1)=C^{\#}W_{/ [R(B)]}C=C^{\#}W(I-Q)C.$$
\end{dem}

\section{Minimization problems in the indefinite trace space}
In the present section the notion of trace of an operator is extended to the Krein space setting with the aim of applying the results previously obtained to trace-type problems on operators.

We denote by $S_p$ the $p${\emph{-Schatten class}} for $1 \leq p < \infty.$ The reader is referred to \cite{Ringrose, Simon} for further details on $S_p$-operators.

Let $(\HH, \K{ \ }{ \ })$ be a Krein space. If $J$ is a signature operator for $\HH,$ fix the Hilbert space $(\HH, \PI{ \ }{ \ }),$ where $\PI{x }{ y }=\K{Jx}{y}$ for all $x, y \in \HH.$ 
The operator $T$ belongs to the Schatten class $S_p(J)$ if $T \in S_p$ when viewed as acting on the associated Hilbert space $(\HH, \PI{ \ }{ \ }).$ The next lemma shows that if $T \in S_p(J_a)$ for some fundamental decomposition of $\HH$ with signature operator $J_a$ then $T \in S_p(J_b)$ for any other fundamental decomposition of $\HH$ with signature operator $J_b.$ To prove this assertion we will use the following result, see \cite[Theorem 2.1.3]{Ringrose}. 

\begin{thm} \label{thmSp} Let $\HH$ be a Hilbert space, $T\in L(\HH)$ and $1 \leq p < \infty.$ Then $T \in S_p$ if and only if there exists a sequence $\{F_n\}_{n \in \mathbb{N}}$ of operators on $\HH$ such that $F_n$ has finite rank not greater than $n$ and $$\sum_{n \geq 1} \Vert T - F_n \Vert^p < \infty.$$
\end{thm} 

\begin{lema} \label{LemaSpKrein} Let $(\HH, \K{ \ }{ \ })$ be a Krein space with signature operators $J_a$ and $J_b$. Fix the Hilbert spaces $(\HH, \PI{ \ }{ \ }_a)$ and $(\HH, \PI{ \ }{ \ }_b).$ Then  $T \in S_p(J_a)$ if and only if $T \in S_p(J_b)$. 
\end{lema}

\begin{dem} The result is readily obtained by applying Theorem \ref{thmSp} and from the fact that $ \PI{ \ }{ \ }_a$ and $ \PI{ \ }{ \ }_b$ are equivalent.
\end{dem}

On account of the above lemma we just write $S_p$ instead of $S_p(J).$

\begin{Def} Let $(\HH, \K{ \ }{ \ })$ be a separable Krein space with signature operator $J$ and fix the associated Hilbert space $(\HH, \PI{ \ }{ \ }).$ If $T \in S_1$ and $\{e_n : n\in \mathbb{N} \}$ is an orthonormal basis of $(\HH, \PI{ \ }{ \ }),$ then the $J${\emph{-trace}} of $T,$ denoted by $\tr_{J}(T),$ is defined as 
	$$\tr_{J}(T)=\sum_{n=1}^{\infty} \K{Te_n}{e_n}.$$
\end{Def}

Notice that $\tr_J(T)$ equals $\tr(JT)$ in the inner product $\PI{ \ }{ \ } =\K{J \ }{ \ }$ see \cite{Ringrose, Simon}. Whence the $J$-trace of $T$ does not depend on the particular choice of the orthonormal basis (see \cite[Lemma 2.2.1]{Ringrose}).

The next lemma gathers the basic properties of the $J$-trace. By using the definition of $\tr_J$ and the properties of the trace of an operator in a Hilbert space the proof is straightforward.

\begin{lema} \label{Proptr} Let $(\HH, \K{ \ }{ \ })$ be a Krein space with signature operator $J$ and fix the associated Hilbert space $(\HH, \PI{ \ }{ \ }).$  Let $T, S \in S_1$ and $\alpha, \beta \in \mathbb{C},$ then
	\begin{enumerate}
		\item [i)] $\tr_{J}(\alpha T + \beta S)= \alpha  \ \tr_{J}(T) + \beta \ \tr_{J}(S),$
		\item [ii)] $\tr_{J}(T^{\#})=\ol{\tr_{J}(T)},$
		\item [iii)] $\tr_{J} (T)= \tr(JT),$ where the trace is calculated with respect to the inner product $\PI{ \ }{ \ } =\K{J \ }{ \ },$
		\item [iv)] $\tr_{J}(TS) = \tr_{J}(JSJT)=\tr_{J}(SJTJ),$
		\item [v)] $\vert \tr_{J}(T) \vert \leq  \Vert T \Vert_{1}.$
	\end{enumerate}
\end{lema}

The next example shows that the $J$-trace depends on the signature operator $J.$

\begin{example} Consider $\mathbb{C}^2$ with the indefinite metric $\K{(x_1, x_2)}{(y_1, y_2)}=x_1\ol{y_1} - x_2 \ol{y_2}.$  
	
	Then $(\mathbb{C}^2, \K{ \ } { \ })$ is a Krein space with fundamental decompositions: $\mathbb{C}^2= span \{ (1,0 )\} \ [\dotplus] \ span \{ (0,1)\}$ and $ \mathbb{C}^2= span \{ (2,1 )\} \ [\dotplus] \ span \{ (1,2)\}.$ Let $J_a$ and $J_b$ be the corresponding signature operators. 
	Observe that $\{ (1,0), (0,1) \}$ is an orthonormal basis in $(\mathbb{C}^2, \K{J_a \ } { \ })$ and $\left\{ \frac{1}{\sqrt{3}}(2,1), \frac{1}{\sqrt{3}}(1,2) \right\}$ is an orthonormal basis in $(\mathbb{C}^2, \K{J_b \ } { \ }).$
	Set $T: \mathbb{C}^2 \ra  \mathbb{C}^2,$ $T(x_1,x_2) :=(x_1+x_2,0).$ A straightforward computation gives $\tr_{J_a}(T)=1 \not = 3= \tr_{J_b}(T).$
	
\end{example}

\begin{lema} \label{propFund} Let $(\HH, \K{ \ }{ \ })$ be a separable Krein space with signature operators $J_a$ and $J_b$. Fix the Hilbert spaces $(\HH, \PI{ \ }{ \ }_a)$ and $(\HH, \PI{ \ }{ \ }_b).$ If $T \in S_1$ then
	$$\tr_{J_b}(T)=\tr_{J_a}(J_b T J_a).$$
\end{lema}

\begin{dem} We use the notation $\tr_{\PI{ \ }{ \ }}$ when we want to highlight the inner product  on which the trace is calculated. Let $\alpha = J_a J_b.$ Then $\alpha$ is an invertible  operator on $\HH$ such that, for every $x, y \in \HH$,
	$$\PI{\alpha x}{y}_a=\PI{J_aJ_bx}{y}_a=\K{J_bx}{y}=\PI{x}{y}_b.$$
	In particular, $\PI{\alpha x}{x}_a \geq 0$ for every $x \in \HH.$ 
	
	Let $\{ e_n : n \in \mathbb{N} \}$ be an orthonormal basis in $(\HH, \PI{ \ }{ \ }_b).$ Then
	$$\delta_{ij}=\PI{e_i}{e_j}_b=\PI{\alpha e_i}{e_j}_a=\PI{\alpha^{1/2}e_i}{\alpha^{1/2}e_j}_a.$$ Hence, $\{ \alpha^{1/2} e_n : n \in \mathbb{N} \}$ is an orthonormal basis in $(\HH, \PI{ \ }{ \ }_a).$ 
	
	\noindent So, if $T \in S_1$ then \\
	$\tr_{J_b}(T)= \tr_{\PI{ \ }{ \ }_b}(J_b T)=\tr_{\PI{ \ }{ \ }_b}(T J_b)=\sum_{ n \geq 1} \PI{TJ_be_n}{e_n}_b=\\=\sum_{ n \geq 1} \PI{\alpha TJ_a \alpha e_n}{e_n}_a =\sum_{ n \geq 1} \PI{ ( \alpha^{1/2} TJ_a \alpha^{1/2})\alpha^{1/2} e_n}{\alpha^{1/2} e_n}_a= \\ =\tr_{\PI{ \ }{ \ }_a}(\alpha^{1/2} T J_a \alpha^{1/2}) = \tr_{\PI{ \ }{ \ }_a}(\alpha T J_a)= \tr_{\PI{ \ }{ \ }_a}(J_a J_b T J_a)=\\
	= \tr_{J_a}(J_b T J_a).$ 
\end{dem}

\subsection*{Fr\'echet derivative of the $J$-trace}

Let  $(\mc{E}, \NC{\cdot})$ be a Banach space and $\mathcal U \subseteq \mc{E}$ be an open set. We recall that a function 
$f: \mc{E} \rightarrow \mathbb{R}$ is said to be {\emph{Fr\'echet differentiable}} at $X_0 \in \mathcal U$ if there exists $Df(X_0): \mc{E} \ra \mathbb{R}$ a bounded linear functional such that
$$\lim\limits_{Y\rightarrow 0} \frac{|f(X_0+Y)-f(X_0) - Df(X_0)(Y)|}{\Vert Y \Vert}=0.$$ 
If $f$ is Fr\'echet differentiable at every $X_0 \in \mc{E}$, $f$ is called Fr\'echet differentiable  on $\mc{E}$  and the function $Df$ which assigns to every point $X_0 \in \mc{E}$ the derivative $Df(X_0),$ is called the Fr\'echet derivative of the function $f.$ If, in addition, the derivative $Df$ is continuous, $f$ is  said to be a {\emph{class $\mc{C}^1$-function}}, in symbols, $f \in \mc{C}^1(\mc{E}, \mathbb{R}).$

%On the other hand, let $t \in \mathbb{R},$ the G\^ateaux derivative of $f$ in  $x \in \mathcal{E}$ in the direction $y \in \mathcal{E}$ eis defined by  
%$$ Df(x,y)=lim_{t \rightarrow 0} \frac{f(x+t y)-f(x)}{t}.$$
\bigskip
Let $W \in L(\HH)^s,$ $B \in CR(\HH)$ and $C\in L(\HH).$ Recall that
$F(X)=(BX-C)^{\#}W(BX-C)$ and consider $f_J: L(\HH) \ra \mathbb{R}$ defined by 
$$f_J(X):=\tr_J(F(X)).$$   	

In the following lemma we give the formula for the Fr\'echet derivative of $f_J(X),$ see \cite{Gold1} for the finite-dimensional case.
\begin{lema} \label{LemaDiff} Let $(\HH, \K{ \ }{ \ })$ be a Krein space with signature operator $J.$ Fix the associated Hilbert space $(\HH, \PI{ \ }{ \ }).$ Let $W \in S_1,$ $B \in CR(\HH)$ and $C \in L(\HH).$
	Then $f_J$ is Fr\'echet differentiable on $L(\HH)$ and $$Df_J(X)(Y)=2 \ Re \ \tr_{J}(Y^{\#}B^{\#}W(BX-C)).$$ Moreover, $f_J \in \mc{C}^1({L(\HH), \mathbb{R}}).$
\end{lema}

\begin{dem} For all  $X, Y \in L(\HH),$ $$f_J(X+Y)=f_J(X)+ 2 Re \ \tr_{J} ((BY)^{\#}W(BX-C)) + \tr_{J} ((BY)^{\#}W(BY)).$$
	Then
	$$\frac{\vert f_J(X+Y)-f_J(X) -  2 Re \ \tr_{J} ((BY)^{\#}W(BX-C))  \vert }{\Vert Y \Vert}  =$$ $$=\frac{\vert \tr_{J} ((BY)^{\#}W(BY)) \vert}{\Vert Y \Vert} \leq \frac{\Vert BY \Vert^2 \Vert W \Vert_1}{\Vert Y \Vert}\leq \Vert B \Vert^2 \Vert W \Vert_{1} \Vert Y \Vert$$ (see Lemma \ref{Proptr}). Hence $f_J$ is Fr\'echet differentiable on $L(\HH)$ and $$Df_J(X)(Y)=2 Re \ \tr_{J} ((BY)^{\#}W(BX-C)).$$ 
	
	Finally, since 
	\begin{align*}
	\vert Df_J(X_1)(Y) - Df_J(X_2)(Y) \vert &= 2 \vert Re \ \tr_{J} ((BY)^{\#}W(B(X_1-X_2))\vert\\
	&\leq 2 \Vert B \Vert^2 \Vert Y \Vert \Vert W \Vert_{1} \Vert X_1-X_2\Vert.
	\end{align*}
	(once again by Lemma \ref{Proptr}), it follows that $f_J \in \mc{C}^1({L(\HH), \mathbb{R}}).$
\end{dem}

\vspace{0,3cm} 
In this section we deal with the following problems: let $(\HH, \K{ \ }{ \ })$ be a Krein space with signature operator $J.$ Fix the associated Hilbert space $(\HH, \PI{ \ }{ \ }).$ Given $B \in CR(\HH),$ $C \in L(\HH)$ and $W \in S_1\cap L(\HH)^s,$ we analyze whether there exists the
\begin{equation}
\underset{X \in L(\HH)}{\min} \tr_{J}((BX-C)^{\#}W(BX-C)) \label{eq71}
\end{equation}
and the corresponding maximum.

Finally, if $B$ is represented as in \eqref{DescomposicionB} and $F_J(X,Y)=(B_+X+B_-Y-C)^{\#}W(B_+X+B_-Y-C),$ we also analyze the existence of
\begin{equation}
\underset{Y \in L(\HH)}{\max} \left(  \underset{X \in L(\HH)}{\min}  \tr_{J}(F_J(X,Y))\right). \label{eq772}
\end{equation}
%
%\vspace{0,3cm} 
It follows from the last lemma that, if $f_J : L(\HH) \times L(\HH) \ra \mathbb{R}$ is given by	
\begin{equation} \label{GJ}
f_J(X,Y):= \tr_J(F_J(X,Y)),
\end{equation}
then $f_J \in \mc{C}^1(L(\HH) \times L(\HH), \mathbb{R})$ and the partial derivatives of $f_J$ in every $(X_0,Y_0) \in L(\HH) \times L(\HH)$ are 
$$D_{X} f_J(X_0,Y_0)(H)=2 Re \ \tr_{J} ((B_+H)^{\#}W(B_+X_0+B_-Y_0-C)), $$
$$D_{Y} f_J(X_0,Y_0)(K)=2 Re \ \tr_{J} ((B_-K)^{\#}W(B_+X_0+B_-Y_0-C)),$$
for all $H, K \in L(\HH).$

%\vspace{0,3cm} 

\begin{thm} \label{thmtrJ} Let $W \in L(\HH)^s$ such that  $W \in S_1,$ $B \in CR(\HH)$ and $C \in L(\HH).$ The following assertions hold: 
	\begin{enumerate}
		\item Assume that $R(B)$ is $W$-nonnegative. Then, $X_0 \in L(\HH)$ realizes  \eqref{eq71} for any signature operator $J$ if and only if  $X_0$ is a $W$-ImS of the equation $BX-C=0.$ 
		\item Let $B$ be represented as in \eqref{DescomposicionB} for some signature operator $J.$ Then, the min-max in \eqref{eq772} exists for every $C \in L(\HH)$ if and only if $W$ is $R(B)$-complementable. In this case,
		$$\underset{Y \in L(\HH)}{\max} \left(  \underset{X \in L(\HH)}{\min}  \tr_{J}(F_J(X,Y))\right)=\tr_J(C^{\#}W_{/[ R(B)]}C).$$
		The operator $Z \in L(\HH)$ realizes \eqref{eq772} if and only if  $Z$ is a $W$-ImMS of $BX-C=0.$
	\end{enumerate}
\end{thm}

\begin{dem} Let  $J$ be a signature operator of $\HH$ and fix the associated Hilbert space $(\HH, \PI{ \ }{ \ }).$  Suppose that $X_0$ is  a solution of Problem \eqref{eq71}. If $f_J$ is as in Lemma \ref{LemaDiff} then $X_0$ is a global minimum of $f_J$. Since $f_J$ is a  $\mc{C}^1$-function, $X_0$ is a critical point of $f_J(X);$ i.e., for every $Y \in L(\HH),$ $Df_J(X_0)(Y)=0$ or equivalently,
	$$0=2 Re \ \tr_{J} ((BY)^{\#}W(BX_0-C))=2 Re \ \tr (J(BY)^{\#}W(BX_0-C)).$$
	Thus, considering a suitable $Y,$ it follows that $$B^{\#}W(BX_0-C)=0.$$ 
	So, by Corollary \ref{cormin2}, $X_0$ is a $W$-ImS of $BX-C=0.$ 
	
	As for the converse, suppose that $X_0$ is a $W$-ImS  of $BX-C=0.$ Let $\{e_n : n\in \mathbb{N} \}$ be any orthonormal basis in $(\HH, \PI{ \ }{ \ }).$ Then
	$$\K{W(BX_0-C)e_n}{(BX_0-C)e_n} \leq \K{W(BX-C)e_n}{(BX-C)e_n}$$ 
	$\mbox{ for every } n \in \mathbb{N} \mbox{ and every } X \in L(\HH).$ Therefore
	$$\tr_{J}(F(X_0)) \leq \tr_{J}(F(X))$$
	for every $X \in L(\HH).$ Hence $X_0$ is a solution of Problem \eqref{eq71} and the proof of the item $1$ is complete.

	As for the item $2,$ suppose that $W$ is $R(B)$-complementable and $Z'$ is a solution of $B^{\#}W(BX-I)=0.$ Then, for any $C\in L(\HH),$ $Z=Z'C$ is a solution of $B^{\#}W(BX-C)=0$ and, by Theorem \ref{TeominmaxW}, $Z$ is a $W$-ImMS of $BX-C=0,$ i.e., 
	$$(BZ-C)^{\#}W(BZ-C)=\underset{Y \in L(\HH) }{\max}  \left( \underset{X \in L(\HH)}{\min} F_J(X,Y)\right).$$ 
	Let  $\{e_n : n\in \mathbb{N} \}$ be any orthonormal basis in $(\HH, \PI{ \ }{ \ }). $ Then, for every $n \in \mathbb{N}$ and any $X, Y \in L(\HH),$ 
	$$\K{(B_+Z+B_-Y-C)^{\#}W(B_+Z+B_-Y-C)e_n}{e_n} \leq$$ $$\leq \K{(B_+Z+B_-Z-C)^{\#}W(B_+Z+B_-Z-C)e_n}{e_n} $$
	$$\leq \K{(B_+X+B_-Z-C)^{\#}W(B_+X+B_-Z-C)e_n}{e_n}.$$
	Therefore
	\begin{align*}
	\tr_{J}(F_J(Z,Z))=\tr_J(F(Z))&=\underset{Y \in L(\HH)}{\max} \left(  \underset{X \in L(\HH)}{\min}  \tr_{J}(F_J(X,Y))\right)=\\
	&=\tr_J(C^{\#}W_{/[ R(B)]}C),
	\end{align*}
	where we used Corollary \ref{CorminmaxregW}. Hence $Z$ is a solution of Problem \eqref{eq772}. 
	
	Conversely, if  $Z \in L(\HH)$ is  a solution of Problem \eqref{eq772} for any $C\in L(\HH),$  then 
	$$f_J(Z,Y) \leq f_J(Z,Z) \leq f_J(X,Z) \mbox{ for every } X,\ Y \in L(\HH),$$ where $f_J$ is as in \eqref{GJ}. Hence, $Z$ is a global minimum of $f_J(X,Z)$ and $Z$ is a global maximum of $f_J(Z,Y).$ Therefore, for every  $H, K \in L(\HH),$ $$D_{X} f_J(Z,Z)(H)=D_{Y} f_J(Z,Z)(K)=0$$ or equivalently,
	
	\noindent $Re \ \tr_{J} ((B_+H)^{\#}W(B_+Z+B_-Z-C))=Re \ \tr_{J} ((B_-K)^{\#}W(B_+Z+B_-Z-C))=0.$ 
	Then, considering suitable $H,K,$ it follows that $$B_+^{\#}W(B_+Z+B_-Z-C)=B_-^{\#}W(B_+Z+B_-Z-C)=0.$$
	Thus 
	$$B^{\#}W(BZ-C)=0$$ and, by Theorem \ref{TeominmaxW} once again,  $Z$ is a $W$-ImMS of $BX-C=0.$ 
\end{dem}

The following theorem synthesizes the results of the last two sections.

\begin{thm} Let $W \in L(\HH)^s$ such that  $W \in S_1$ and $B \in CR(\HH).$ Then the following statements are equivalent:
	\begin{itemize}
		\item [i)] there exists a $W$-ImMS of  $BX-C=0$ for every $C \in L(\HH),$
		\item [ii)] the $\underset{Y \in L(\HH)}{\max} \left(  \underset{X \in L(\HH)}{\min}  \tr_{J}(F_J(X,Y))\right)$ is attained, for every $C \in L(\HH),$ 
		\item [iii)]  $W$ is $R(B)$-complementable,
		\item [iv)] the equation $B^{\#}W(BX-C)=0$ admits a solution for every $C \in L(\HH).$ 
	\end{itemize}
	In this case, $$\underset{Y \in L(\HH)}{\max} \left( \underset{X \in L(\HH)}{\min} F_J(X,Y)\right)=C^{\#}W_{/ [R(B)]}C$$
	and $$\underset{Y \in L(\HH)}{\max} \left(  \underset{X \in L(\HH)}{\min}  \tr_{J}(F_J(X,Y))\right)=\tr_J(C^{\#}W_{/[ R(B)]}C).$$
	
	Moreover, $Z$ is a $W$-ImMS of $BX-C=0$ and the min-max in $ii)$ is attained in $Z$ if and only $Z=Z_1+Z_2,$ where $B^{\#}W(BZ_1-C)=0$ and $(BZ_2)^{\#}WBZ_2=0.$
	
\end{thm}

\subsection*{\textbf{Final remark: the $\mathbf{J}$-$\mathbf{S_2}$ space}}
Let $(\HH, \K{ \ }{ \ })$ be a Krein space with signature operator $J.$ Fix the associated Hilbert space $(\HH, \PI{ \ }{ \ })$ and set
$$\K{S}{T}_J:=\tr_{J}(T^{\#}S), \quad S,T \in S_2.$$ 
It can be readily seen that $\K{ \ } { \ }_J$ is an indefinite inner product on $S_2.$ Moreover, $(S_2, \K{ \ } { \ }_J)$ is a Krein space and 
$$\tr_{J}(T^{\#}T)=\Vert P_+ T \Vert_2^2-\Vert P_- T \Vert_2^2,$$ where $P_{\pm}=\frac{I\pm J}{2}.$

\section*{Acknowledgements}
Maximiliano Contino and Alejandra Maestripieri were supported by CONICET PIP 0168.  The work of Stefania Marcantognini was done during her stay at  the Instituto Argentino de Matem\'atica with an appointment funded by the CONICET. She is greatly grateful to the institute for its hospitality and to the CONICET for financing her post.


\section*{References}
\begin{thebibliography}{999}
	
	\bibitem{Shorted2} Anderson~W.N., Trapp~G.E., {\it Shorted Operators II}, SIAM J. Appl. Math., 28 (1975), 60-71.
	
	%\bibitem{AndoSchur} Ando~T., {\it Generalized Schur complements}, Linear Algebra Appl., 27 (1979), 173-186.
	
	\bibitem{AndoLibro} Ando~T., {\it Linear operators on Krein spaces}, Hokkaido University, Sapporo, Japan (1979).
	
	\bibitem{AntCorSto06} Antezana~J., Corach~G., Stojanoff~D., {\it Bilateral shorted operators and parallel sums}, Linear Algebra Appl., 414 (2006), 570-588.
	
	\bibitem{Azizov} Azizov~T.Y., Iokhvidov~I.S., {\it Linear operators in spaces with and indefinite metric}, John Wiley and Sons, 1989.
	
	%\bibitem{Baidiuk} Baidiuk~D., {\it Completion and extension of operators in Krei?n spaces}, J. Math. Sci., 224 (2017), 493-508.
	
	%\bibitem{Baidiuk2} Baidiuk~D., Hassi~S., {\it Completion, extension, factorization, and lifting of operators}, Math. Ann., 364 (2016), 1415-1450.
	
	\bibitem{Bognar} Bogn\'ar~J., {\it Indefinite inner product spaces}, Springer, Berlin (1974).
	
	%\bibitem{BognarKramli} Bogn\'ar~J., Kr\'amli~A., {\it Operator of the form $C^*C$ in indefinite product spaces}, Acta Sci. Math. (Szeged), 29 (1968), 19-29.
	
	%\bibitem{Carlson} Carlson~D., Haynsworth~E.V., {\it Complementable and almost definite matrices}, Linear Algebra Appl., 52 (1983), 157-176.
	
	\bibitem{CMSSzeged} Corach~G., Maestripieri~A., Stojanoff~D., {\it Oblique projections and Schur complements}, Acta Sci. Math. (Szeged), 67 (2001), 337-256.
	
	\bibitem{Contino} Contino~M., Giribet~J.I., Maestripieri~A., {\it Weighted Procrustes problems}, J. Math. Anal. Appl., 445 (2017), 443-458.
	
	%\bibitem{Contino2} Contino~M., Giribet~J.I., Maestripieri~A., {\it Weighted least square solutions of the equation $AXB-C=0$}, Linear Algebra Appl. 518 (2017), 177-197.
	
	\bibitem{Contino3} Contino~M., Maestripieri~A., Marcantognini~S., {\it Operator least squares problems and Moore-Penrose inverse in Krein Spaces}, Integr. Equat. Oper. Th., 90 (2018), 32. 
	
	\bibitem{Contino4} Contino~M., Maestripieri~A., Marcantognini~S., {\it Schur complements of selfadjoint Krein space operators}, (2018), arXiv:1809.01695.
	
	\bibitem{WGI} Corach~G., Fongi~G., Maestripieri~A., {\it Weighted projections into closed subspaces}, Studia Mathematica, 216 (2013), 131-148.
	
	\bibitem {Douglas} Douglas~R.G., {\it On majorization, factorization and range inclusion of operators in Hilbert space}, Proc. Amer. Math. Soc., 17 (1966), 413-416.
	
	\bibitem{DR} Dritschel~M.A., Rovnyak~J., {\it Extension theorems for contraction operators on Krein spaces}, Operator Theory: Adv. Appl., 47 (1990), 221-305. 
	
	\bibitem{DR1} Dritschel~M.A., Rovnyak~J., {\it Operators on indefinite inner product spaces}, Lectures on operator theory and its applications, 3 (1996), 141-232.
	
	\bibitem{Nashed}  Engl~H.W., Nashed~M.Z., {\it New extremal characterizations of generalized inverses of linear operators}, J. Math. Anal. Appl., 82 (1981), 566-586.
	
	
	%\bibitem{GiribetKreinShorted}  Giribet J.I., Maestripieri A., Mart\'inez Per\'ia F., {\it Shorting selfadjoint operators in Hilbert spaces}, Linear Algebra Appl., 428 (2008), 1899-1911.
	
	%\bibitem{GiribetIndefinite} Giribet~J.I., Maestripieri~A., Mart\'inez Per\'ia~F., {\it Indefinite least-squares problems and pseudo-regularity},  J. Math. Anal. Appl., 430 (2015), 895-908.
	
	\bibitem{GiribetKrein} Giribet~J.I., Maestripieri~A., Mart\'inez Per\'ia~F., {\it A geometrical approach to indefinite least squares problems}, Acts Appl. Math, 111 (2010), 65-81.
	
	\bibitem{Gold1} Goldstein ~G.R., Goldstein~J.A., {\it The best generalized inverse}, J. Math. Anal. Appl., 252 (2000), 91-101.
	%\bibitem{Hassi} Hassi S., Nordstr\"om K., {\it On projections in a space with an indefinite metric}, Linear Algebra and its Applications, 208-209 (1974), 401-417.
	
	\bibitem{HassibipartI} Hassibi~B.,  Sayed~A.H., Kailath~T., {\it Linear estimation in Krein spaces - part I: theory}, IEEE Trans. Automat. Control 41 (1996) 18-33.
	
	\bibitem{HassibipartII} Hassibi~B.,  Sayed~A.H., Kailath~T., {\it Linear estimation in Krein spaces - part II: applications}, IEEE Trans. Automat. Control 41 (1996) 33-49.
	
	
	%\bibitem{HassibiIII} Hassibi B.,  Sayed A. H., Kailath T., {\it Indefinite-Quadratic Estimation and Control. A Unified Approach to  $\HH^2$ and $\HH^\infty$ Theories}, Studies in Applied and Numerical Mathematics, 1999.
	
	
	\bibitem{Ulric} Kintzel~U., {\it Procrustes problems in finite dimensional indefinite scalar product spaces}, Linear Algebra Appl., 402 (2005), 1-28.
	
	\bibitem {Krein} Krein~M.G.,  {\it The theory of self-adjoint extensions of semibounded Hermitian operators and its applications}, Mat. Sb. (N.S.), 20 (62) (1947), 431-495.
	
	%\bibitem{KreinSzeged} Maestripieri~A., Mart\'inez Per\'ia~F., {\it Decomposition of selfadjoint projections in Krein spaces}, Acta Sci. Math.(Szeged), 72 (2006), 611-638.
	
	%\bibitem{SchurKrein} Maestripieri~A., Mart\'inez Per\'ia~F., {\it Schur complements in Krein spaces}, Integr. Equ. Oper. Theory, 59 (2007), 207-221.
	
	%\bibitem{XavierMary} Mary~X., {\it Moore-Penrose inverse in Krein spaces}, Integr. Equ. Oper. Theory, 60 (2008), 419-433.
	
	%\bibitem{Massey} Massey~P., Stojanoff~D., {\it Generalized Schur complements and $P$-complementable operators}, Linear Algebra Appl., 393 (2004), 299-318.
	
	\bibitem{Mitra} Mitra~S.K., Rao~C.R., {\it Projections under seminorms and generalized Moore Penrose inverses and operator ranges}, Linear Algebra Appl., 9 (1974), 155-167.
	
	%\bibitem{Pedersen} Pedersen~G.K., {\it Some operator monotone functions}, Proc. Amer. Math. Soc., 36 (1972), 309-310.
	
	%\bibitem{Pekarev} Pekarev~E.L., {\it Shorts of operators and some extremal problems}, Acta Sci. Math. (Szeged), 56 (1992), 147-163.
	
	\bibitem{Ringrose} Ringrose~J.R., {\it Compact non-self-adjoint operators}, Van Nostrand Reinhold Co., 1971.
	
	\bibitem{Hassibietal} Sayed~A.H., Hassibi~B., Kailath~T., {\it Inertia conditions for the minimization of quadratic forms in indefinite metric spaces}, Operator Theory: Adv Appl., 87 (1996), 309-347.
	
	\bibitem{Simon} Simon~B., {\it Trace Ideals and their applications}, London Mathematical Society Lecture Note Series, vol. 35, Cambridge University Press, Cambridge, 1979.
	
	%\bibitem{Pseudo} Maestripieri, A., Per\'ia, F. M. { \it Normal projections in Krein spaces}, Integral Equations Operator Theory, 76 (2013), 357-380.
	
	
	%\bibitem{Nam80} Nambooripad S. K. S. , {\it The natural partial order on a regular semigroup}, Proceedings of the Edinburgh Mathematical Society (Series 2) 23 03 (1980), 249-260.
	
	%\bibitem{Mit86}  Mitra S.K., {\it The minus partial order and the shorted matrix}, Linear Algebra Appl. 83 (1986), 1-27.
	
	%\bibitem{Aiken}  Aiken J.G., Erdos J.A., Goldstein J.A.,  {\it Unitary Approximation of positive operators}, Illinois J. Math, 20 (1980),  61-72.
	
	%\bibitem{AntCorSto06} Antezana J.,    Corach G.,   Stojanoff D. {\it Bilateral shorted operators and parallel sums}, Linear Algebra Appl. 414 2 (2006), 570-588.
	
	%\bibitem{Aiken2} Aiken J.G., Erdos J.A., Goldstein J.A., {\it On L\"owdin Orthogonalization}, International Journal of Quantum Chemistry, 18 (1980), 1101-1108.
	%
	%\bibitem{Bounkhel} Bounkhel, M., {\it On minimizing the norm of linear maps in $Cp$-classes}, Applied Sciences, 8 (2006), 40-47.
	%
	
	%\bibitem{Shorted1} Anderson W.N., {\it Shorted Operators}, SIAM J. Appl. Math, 20 (1971), 520-525.
	
	
	%\bibitem{Kitta}  Alizadeh R., Asadi M. B., {\it An extension of Ky Fan's dominance theorem}, Banach J. Math. Anal., 6 (2012), 139-146.
	
	%\bibitem{Arias} Arias M.L., Gonzalez M.C., {\it Positive solutions to operator equations $AXB=C$}, Linear Algebra and its Applications, 433 (2010), 1194-1202.
	
	%\bibitem{Changsen} Changsen Y., {\it On the critical points of the map 
	%$F_p: X \ra \Vert AXB-C \Vert_{p}^p$}, Applied Mathematics and Mechanics, 21 (2000), 485-488.
	
	%\bibitem{SchurSzeged} Corach G., Maestripieri A., Stojanoff D., {\it Generalized Schur complements and oblique projections}, Linear Algebra Appl., 341 (2002), 259-272.
	
	
	%\bibitem{Shorted3} Corach G., Maestripieri A., Stojanoff D., {\it Generalized orthogonal projections and shorted operators}, Margarita Mathem\'atica, Departamento de Matem\'aticas y Computaci\'on, Universidad de La Rioja, (2001), 607-625.
	%
	%\bibitem{Spline} Corach G., Maestripieri A., Stojanoff D., {\it Oblique projections and abstract splines}, Journal of Approximation Theory, 117 (2002), 189-206.
	
	%\bibitem {Minus} Djiki\'c, M.S., Fongi G., Maestripieri A., {\it Minus order and range additivity}.
	
	%\bibitem {Drago} Dragoljuf J. Kecki\'c, {\it Orthogonality in $C_1$ and $C_{\infty}$ spaces and normal derivations}, J. Operator Theory, 51 (2004),  89-104.
	
	%\bibitem {Eldar} Eldar, Y. C., Werther, T., {\it General framework for consistent sampling in Hilbert spaces}, International Journal of Wavelets, Multiresolution and Information Processing, 3 (2005), 497-509.
	
	%\bibitem{Har80} Hartwig R.E., {\it How to partially order regular elements}, Math. Jpn 25 (1980) 1-13.
	
	
	
	%\bibitem{Nashed2} Nashed, M. Z. {\it Inner, outer, and generalized inverses in Banach and Hilbert spaces}, Numer. Funct. Anal. Optim. 9 (1987), 261-325.
	
	%\bibitem{GiribetKreinShorted}  Giribet J.I, Maestripieri A., Mart\'inez Per\'ia F., {\it Shorting selfadjoint operators in Hilbert spaces}, Linear Algebra Appl., 428 (2008), 1899-1911.
	
	%\bibitem{Mah1} Maher P. J.,  {\it Some matrix approximation problems arising from quantum chemistry}, Proc. Indian Nat. Sci. Acad., Part A 64 (1998), 715-723.
	%\bibitem{Mah3} Maher P. J.,  {\it Some operator inequalities concerning generalized inverses}, Illinois Journal of Mathematics (1990), Vol 34, No. 3.
	%
	%\bibitem{Mah4} Maher P. J., {\it Some norm inequalities concerning generalized inverses, 2}, Linear Algebra and its Applications, 420 (2007), 517-525.
	%
	%\bibitem{Mecheri} Mecheri S., Bounkhel M., {\it Global minimum and orthogonality in $C_1$-classes}, J. Math. Anal. Appl., 287 1 (2003), 51-60. 
	%
	%\bibitem{Mecheri2} Mecheri S., {\it Global minimum and orthogonality in $C_p$-classes}, Math. Nachr., 280 (2007), 794-801.
	
	%\bibitem{Sem10} $\check{\textrm{S}}$emrl P., {\it Automorphisms of $B (\HH)$ with respect to minus partial order}, Journal of Mathematical Analysis and Applications 369 1 (2010), 205-213.
	%%\bibitem{Unser} Unser, M., {\it Sampling-50 years after Shannon}, Proceedings of the IEEE, 88 (2000), 569-587.
	
\end{thebibliography}
\end{document}